\UseRawInputEncoding
\documentclass[12pt]{article}
\usepackage[centertags]{amsmath}
\usepackage{amsfonts}
\usepackage{amssymb}
\usepackage{latexsym}
\usepackage{amsthm}
\usepackage{newlfont}
\usepackage{graphicx}
\usepackage{listings}
\usepackage{booktabs}
\usepackage{abstract}
\usepackage{enumerate}
\usepackage{xcolor}
\RequirePackage{srcltx}
\date{}

\newlength{\defbaselineskip}
\newcommand{\setlinespacing}[1]%
           {\setlength{\baselineskip}{#1 \defbaselineskip}}

\newcommand{\actaqed}{\hfill $\actabox$}
{\medskip\noindent \textit{Proof of #1. }}%
{\actaqed \medskip}

\def\cD{{\mathcal D}}

\def\cS{{\mathcal S}}

\def\bfe{\mathbf e}

\def\bv{\mathbf v}

 \def \<{\langle}
\def\>{\rangle}

\def \e{\varepsilon}

\def \de{\delta}

\def \ff{\varphi}

\def\bt{\beta}

\def \sp{\operatorname{span}}
\def \csp{\overline{\operatorname{span}}}

\def\bt{\beta}

\newtheorem{Theorem}{Theorem}[section]
\newtheorem{Lemma}{Lemma}[section]

\newtheorem{Corollary}{Corollary}[section]
\numberwithin{equation}{section}

\newcommand{\be}{\begin{equation}}
\newcommand{\ee}{\end{equation}}

\begin{document}

\title{On stability of Weak Greedy Algorithm in the presence of noise} 

\author{V.N. Temlyakov \footnote{
		This work was supported by the Russian Science Foundation under grant no. 23-71-30001, https://rscf.ru/project/23-71-30001/, and performed at Lomonosov Moscow State University.
  }}

\newcommand{\Addresses}{{
  \bigskip
  \footnotesize

  V.N. Temlyakov, \textsc{University of South Carolina, USA,\\ Steklov Mathematical Institute of Russian Academy of Sciences, Russia;\\ Lomonosov Moscow State University, Russia; \\ Moscow Center of Fundamental and Applied Mathematics, Russia;\\  Saint Petersburg State University,  Saint Petersburg, Russia. \\ 
  \\
E-mail:} \texttt{temlyakovv@gmail.com}

}}

\maketitle

\begin{abstract}{This paper is devoted to the theoretical study of the efficiency, namely, stability of some greedy algorithms. In the greedy approximation theory researchers are mostly interested in the following   two important properties of an algorithm -- convergence and rate of convergence. 
In this paper we present some results on one more important property of an algorithm -- stability. Stability   means that small perturbations do not result in a large change in the outcome of the algorithm. In this paper we discuss one kind of perturbations -- noisy data.
 	}
\end{abstract}


\section{Introduction}
\label{In}

This paper is devoted to the theoretical study of the efficiency, namely, stability of some greedy algorithms. Greedy algorithms provide sparse approximations with respect to a given 
dictionary (system). Sparse approximation is important in many applications because of concise form of an approximant and 
good accuracy guarantees. The theory of compressed sensing, which proved to be very useful in the image processing and data sciences, is based on the concept of sparsity. A fundamental issue of sparse approximation is the problem of construction of efficient algorithms, which provide good approximation. It turns out that greedy algorithms with respect to dictionaries are very good from this point of view. They are simple in implementation and there are well developed theoretical guarantees of their efficiency. 
The reader can find results on greedy approximation in the books \cite{VTbook}, \cite{VTbookMA}, and a very recent survey \cite{VT211}. Greedy algorithms are very useful in applications; in particular, adaptive methods are used 
in PDE solvers, and sparse approximation is used in image/signal/data processing, as well as in the design of neural networks, and in convex optimization. This fact motivated deep theoretical study of a variety of greedy algorithms. In this paper, we study  the Pure Greedy Algorithm (PGA) and its modifications from the point of view of stability. 

We proceed to a systematic presentation of known and new results. Let $H$ be a real Hilbert space with the inner product $\<\cdot,\cdot\>$ and norm $\|\cdot\|$.  We say that a set of elements (functions) $\cD$ from $H$ is a dictionary (symmetric dictionary) if each $g\in \cD$ has norm   one ($\|g\|= 1$) and $\csp \cD =H$. In addition, we assume for convenience the property of symmetry:
$$
g\in \cD \quad \text{implies} \quad -g \in \cD.
$$
Note that we assume that $H$ is a real Hilbert space and that $\cD$ is symmetric for convenience. The results and arguments below apply to the case of complex Hilbert spaces and any dictionary $\cD$ after a proper straightforward modification. 

We define the Pure Greedy Algorithm (PGA). This algorithm (under another name) was introduced in \cite{FS} for a specific dictionary -- the ridge dictionary. We describe this algorithm for a general dictionary $\cD$. If $f\in H$,
we let $g(f)\in \cD$ be an element from $\cD$ which maximizes $\< f,g\>$. We assume for simplicity that such a maximizer exists; if not, suitable modifications are necessary (see Weak Greedy Algorithm below) in the algorithm that follows. We define:
$$
G(f,\cD):= \<f,g(f)\>g(f)\quad \text{and}\quad R(f,\cD) := f-G(f,\cD).
$$
 
 {\bf Pure Greedy Algorithm (PGA).}   We define $f_0:= f$ and $G_0(f,\cD) := 0$. Then, for each $m\ge 1$, we inductively define
$$
G_m(f,\cD):= G_{m-1}(f,\cD) +G(f_{m-1},\cD),
$$
$$
f_m:= f-G_m(f,\cD) = R(f_{m-1},\cD).
$$

Note that for a given element $f$, the sequence $\{G_m(f,\cD)\}$ may not be unique. 
We now define the most general modification of the PGA, which we study in this paper. 
  Let a sequence $\tau = \{t_k\}_{k=1}^\infty$, $0\le t_k \le 1$, and a parameter $b\in (0,1]$ be given.  We define the Weak Greedy Algorithm with parameter $b$.  

{\bf Weak Greedy Algorithm with parameter $b$ (WGA($\tau,b$)).} We define $f_0:=f_0^{\tau,b}:=f$. Then, for each $m\ge 1$, we inductively define:
\begin{enumerate}
\item[(1)] $\ff_m:=\varphi^{\tau,b}_m \in \cD$ is any satisfying: 
\end{enumerate}
$$
\<f_{m-1},\varphi_m\> \ge t_m \sup_{g\in \cD} \<f_{m-1},g\>;
$$
\begin{enumerate}
\item[(2)]
\end{enumerate}
$$
f_m :=f_m^{\tau,b}:= f_{m-1} -b\<f_{m-1},\varphi_m\>\varphi_m;
$$

\begin{enumerate}
\item[(3)]
\end{enumerate}
$$
G_m(f,\cD):=G^{\tau,b}_m(f,\cD) := b\sum_{j=1}^m \<f_{j-1},\varphi_j\>\varphi_j.
$$
 In the case $t_k = t$, $k=1,2,\dots$, we write $t$ in the notation instead of $\tau$. 

In the case $b=1$ the WGA($\tau,1$) is called the Weak Greedy Algorithm (WGA($\tau$)). 
It was introduced in \cite{VT75}. The WGA($\tau,b$) was introduced in \cite{VT111}. The following Theorem \ref{InT1} was proved in \cite{VT111}. For a real Hilbert space $H$ and a dictionary $\cD\subset H$ define $A_1(\cD)$ to be the closure (in $H$) of the convex hull of the symmetrized dictionary $\cD^\pm := \{\pm g : g\in \cD\}$. 

\begin{Theorem}[{\cite{VT111}}]\label{InT1} Let $\cD$ be an arbitrary dictionary in $H$. Assume $\tau :=\{t_k\}_{k=1}^\infty$, $t_k \in (0,1]$, $k=1,2,\dots$, is a nonincreasing sequence and $b\in(0,1]$. Then, for $f \in A_1(\cD)$, we have:
\be\label{In1}
\|f -G^{\tau,b}_m(f,\cD)\| \le e_m(\tau,b),
\ee
where:
\be\label{In2}
e_m(\tau,b) := \left(1+b(2-b)\sum_{k=1}^m t^2_k\right)^{-\frac{(2-b)t_m}{2(2+(2-b)t_m)}}. 
\ee
\end{Theorem}

In the greedy approximation theory researchers are mostly interested in the following   two important properties of an algorithm -- convergence and rate of convergence. 
In this paper we present some results on one more important property of an algorithm -- stability. Stability can be understood in different ways. Clearly, stability means that small perturbations do not result in a large change in the outcome of the algorithm. In this section we discuss one kind of perturbations -- noisy data. We present some comments in the general setting of greedy approximation in Banach spaces. We only formulate some results for the reader to get a feeling of them and refer the reader for the details to \cite{VT211}, p.48, section Stability. Also, some comments are given in Section \ref{D}. 
Usually, in the greedy algorithms literature the noisy data is understood in the following deterministic way. For a real Banach space $X$ and a dictionary $\cD\subset X$ define $A_1(\cD)$ to be the closure (in $X$) of the convex hull of the symmetrized dictionary $\cD^\pm := \{\pm g : g\in \cD\}$. For each $f\in X$ we associate the following norm
$$
\|f\|_{A_1(\cD)} := \inf \{M:\, f/M\in A_1(\cD)\}.
$$

Take a number $\e\ge 0$ and two elements $f$, $f^\e$ from a Banach space $X$ such that
\be\label{S1}
\|f-f^\e\| \le \e,\quad
f^\e/B \in A_1(\cD),
\ee
with some number $B>0$.
Then we interpret $f$ as a noisy version of our original signal $f^\e$, for which we know that it has some good properties formulated in terms of $A_1(\cD)$. The first results on approximation of noisy data (in the sense of (\ref{S1})) were obtained in \cite{VT115} for the Weak Chebyshev Greedy Algorithm (WCGA) and the Weak Greedy Algorithm with Free Relaxation (WGAFR). The WCGA is the generalization to the case of Banach spaces of the Weak Orthogonal Greedy Algorithm (WOGA) defined in Hilbert spaces (see Section \ref{D} below). 
Later, in \cite{VT165} we proved Theorem \ref{ST1}, which 
covers all algorithms from the collection Weak Biorthogonal Greedy Algorithms (WBGA($\tau$)). Thus, Theorem \ref{ST1} covers the known results from  \cite{VT115} on WCGA and WGAFR and also covers the corresponding results on some other greedy algorithms.

\begin{Theorem}[{\cite{VT165}}]\label{ST1} Let $X$ be a uniformly smooth Banach space with modulus of smoothness $\rho(u)\le \gamma u^q$, $1<q\le 2$. Assume that $f$ and $f^\e$ satisfy (\ref{S1}).
Then, for any algorithm from the collection  WBGA($\tau$) applied to $f$ we have
$$
\|f_m\| \le  \max\left\{2\e,\, C(q,\gamma)(B+\e) \Big(1+\sum_{k=1}^mt_k^p\Big)^{-1/p}\right\},
\quad p:=\frac{q}{q-1},
$$
with $C(q,\gamma)= 4(2\gamma)^{1/q}$.
\end{Theorem}

\section{Stability result for the WGA}
\label{St}

We proceed to the new results of the paper.

\begin{Theorem}\label{T4.2} Let $\cD$ be an arbitrary dictionary in $H$. Assume $\tau :=\{t_k\}_{k=1}^\infty$, $t_k \in (0,1]$, $k=1,2,\dots$, is a nonincreasing sequence and $b\in(0,1]$. Take a number $h \in (0,1)$ and denote
$$
\bt_k := \left(1-\frac{b}{2}\right)ht_k.
$$
Suppose that for some $\e\in (0,1]$  $f$ and $f^\e$ are such that 
\be\label{St1}
\|f-f^\e\|\le \e,\quad f^\e/B \in A_1(\cD),  \quad B>0.
\ee
Then   we have after $m\le \e^{-2}$  iterations of the WGA$(\tau,d)$
\be\label{St4.1}
 \|f^{\tau,b}_m\| \le \max\left(\frac{\e}{1-h},\|f\|^{\frac{1}{1+\bt_m}}(B+1)^{\frac{\bt_1}{1+\bt_m}}e_m(\tau,b,h) \right), 
\ee
where
$$
e_m(\tau,b,h) := \left(\left(\frac{h\|f\|}{B+1}\right)^{-2}+ b(2-b) \sum_{k=1}^m t_k^2\right)^{-\frac{\bt_m}{2(1+\bt_m)}}.
$$
\end{Theorem}
\begin{Corollary}\label{StC1} Consider the special case $\tau =\{t\}$, $t\in (0,1]$. Our assumption (\ref{St1}) implies $\|f\| \le B+\e \le B+1$. 
Denote
$$
\bt:=\bt(t,b,h) := \left(1-\frac{b}{2}\right)ht.
$$
Then, under the conditions of Theorem \ref{T4.2}  we have
\be\label{St4.1}
\|f^{t,b}_m\| \le \max\left(\frac{\e}{1-h},(B+1) \left(b(2-b) m h^2 t^2\right)^{-\frac{\bt}{2(1+\bt)}}\right).  
\ee
\end{Corollary}

Corollary \ref{StC1} shows that after   
\be\label{m}
 m\in [(2\e)^{-2},\e^{-2}]
\ee
iterations we obtain 
\be\label{er1}
\|f^{t,b}_m\| \le C_1\e^{\frac{\bt}{1+\bt}} 
\ee
with a positive constant $C_1$ independent of $\e$. 
Note that Theorem \ref{InT1} provides the following error bound for the noiseless signal 
$f^\e$ after the same number (see (\ref{m})) of iterations
\be\label{er2}
\|(f^\e)^{t,b}_m\| \le C_2 \e^{\frac{\bt((t,b,1))}{1+\bt(t,b,1)}}
\ee
with a positive constant $C_2$ independent of $\e$. 
It is clear that $\bt(t,b,h)$ is close to $\bt(t,b,1)$, when $h$ is close to $1$ and the bound (\ref{er2}) is close to the bound (\ref{er1}). 

\begin{proof}  We introduce some notations: for brevity $f_m:=f^{\tau,b}_m$, $\varphi_m:=\varphi^{\tau,b}_m$ and
$$
a_m := \|f_m\|^2, \quad y_m := 
\<f_{m-1},\varphi_m\>, \quad
m=1,2,\dots,  
$$
and consider the sequence $\{B_n\}$ defined as follows
$$
B_0 := B+1,\quad B_m := B_{m-1} +by_m, \quad m=1,2,\dots .
$$
Along with the sequence 
$$
f_m := f - b\sum_{k=1}^m \<f_{k-1},\ff_{k}\>\ff_k
$$
consider the sequence $(f^\e_0:=f^\e)$
$$
f_m^\e := f^\e - b\sum_{k=1}^m \<f_{k-1}^\e,\ff_{k}\>\ff_k,\quad f_{k}^\e:= f_{k-1}^\e- b\<f_{k-1}^\e,\ff_{k}\>\ff_k, \quad k=1,2,\dots.
$$
Then, using the notation $\de_k:= f_k - f_k^\e $, $k=0,1,\dots$, we obtain
\be\label{St2}
\de_k   = \de_{k-1}  - b\<\de_{k-1},\ff_{k}\>\ff_{k}
\ee
and
\be\label{St2a}
\|\de_k\|^2   = \|\de_{k-1}\|^2  - b(2-b)\<\de_{k-1},\ff_{k}\>^2.
\ee
This and our assumption (\ref{St1}) imply
\be\label{St3}
\|\de_k\| \le \|\de_{k-1}\| \le \dots \le \|\de_0\|\le \e.
\ee
Also (\ref{St2a}) implies
\be\label{St3a}
   b(2-b)\sum_{j=1}^k\<\de_{j-1},\ff_{j}\>^2 \le \|\de_0\|^2 \le \e^2.
\ee
By Lemma 6.10 from \cite{VTbook}, p. 343, we obtain
$$
\sup_{g\in \cD} \<f_{m-1},g\> = \sup_{\phi\in A_1(\cD)}  \<f_{m-1},\phi\> \ge \<f_{m-1},f^\e_{m-1}\>\|f^\e_{m-1}\|_{A_1(\cD)}^{-1}
$$
$$
= \<f_{m-1},f_{m-1} +f^\e_{m-1}-f_{m-1}\>\|f^\e_{m-1}\|_{A_1(\cD)}^{-1},
$$
using (\ref{St3}) we continue
\be\label{St4}
 \ge (\|f_{m-1}\|^2 -\e\|f_{m-1}\|)\|f^\e_{m-1}\|_{A_1(\cD)}^{-1}.
\ee
Thus under assumption $(1-h)\|f_{m-1}\|\ge \e$ we get
\be\label{St5}
\sup_{g\in \cD} \<f_{m-1},g\> \ge \frac{h\|f_{m-1}\|^2}{\|f^\e_{m-1}\|_{A_1(\cD)}}.
\ee
It is clear that for $k\le \e^{-2}$ we have $\|f^\e_{k}\|_{A_1(\cD)} \le  B_{k}$. Indeed, $\|f^\e\|_{A_1(\cD)} \le B$ and for $k\ge 1$ we have
$$
\|f^\e_k\|_{A_1(\cD)}\le \|f^\e\|_{A_1(\cD)}+ b \sum_{j=1}^k|\<f_{j-1},\ff_{j}\>| + b\sum_{j=1}^k|\<\de_{j-1},\ff_{j}\>|
$$
and by (\ref{St3a}) continue
$$
\le B_k-1 + bk^{1/2} \left(\sum_{j=1}^k\<\de_{j-1},\ff_{j}\>^2\right)^{1/2} \le  B_k-1 +k^{1/2}b^{1/2}\e \le B_k. 
$$
Therefore, we obtain
\be\label{St4.2}
\sup_{g\in \cD} \<f_{m-1},g\> \ge  \frac{ha_{m-1}}{B_{m-1}}.
\ee
From here and from the equality  
$$
\|f_m\|^2 = \|f_{m-1}\|^2 -b(2-b)\<f_{m-1},\varphi_m\>^2
$$
we obtain the following relations
\be\label{St4.3}
a_m = a_{m-1} - b(2-b)y_m^2, 
\ee
\be\label{St4.4}
B_m = B_{m-1} +by_m,  
\ee
\be\label{St4.5}
y_m \ge ht_ma_{m-1}/B_{m-1}.  
\ee
From (\ref{St4.3}) and (\ref{St4.5}) we get
$$
a_m \le a_{m-1} (1-b(2-b)(ht_m)^2a_{m-1}B_{m-1}^{-2}).
$$ 
Using that $B_{m-1}\le B_m$ we derive from here
$$
a_mB_m^{-2} \le a_{m-1}B_{m-1}^{-2}(1-b(2-b)(ht_m)^2a_{m-1}B_{m-1}^{-2}).
$$
We now use the following known Lemma \ref{HL1} (see, for instance, \cite{VT211}, Lemma 14.2).

\begin{Lemma}\label{HL1} Let a number $C>0$ and a sequence $\{v_k\}_{k=1}^\infty$, $v_k \ge 0$, $k=1,2,\dots$, be given.
Assume that $\{x_m\}_{m=0}^\infty$
is a sequence of non-negative
 numbers satisfying the inequalities
$$
x_0 \le C, \quad x_{m} \le x_{m-1}(1 - x_{m-1} v_{m}) , \quad m = 1,2, \dots  .
$$
Then we have for each $m$
$$
x_m \le   \left(C^{-1}+ \sum_{k=1}^{m} v_k\right)^{-1} .
$$
\end{Lemma}

By Lemma \ref{HL1}  with  $x_m:=a_mB_m^{-2}$, $C := (\|f\|/B_0)^2$, $v_m:=b(2-b)(ht_m)^2$  we obtain
\be\label{St4.6}
a_mB_m^{-2} \le   \left(C^{-1}+  b(2-b)\sum_{k=1}^m (ht_k)^2\right)^{-1}.  
\ee
The relations (\ref{St4.3}) and (\ref{St4.5})  imply
$$
a_m \le a_{m-1}-b(2-b)y_mht_ma_{m-1}/B_{m-1} 
$$
\be\label{St4.7}
 = a_{m-1}(1-b(2-b)ht_my_m/B_{m-1}).  
\ee
Using the inequality $(1-x)^{1/2} \le 1-x/2$ we get from (\ref{St4.7})
\be\label{St4.8}
a_m^{1/2} \le a_{m-1}^{1/2}(1-b(1-b/2)ht_my_m/B_{m-1})=a_{m-1}^{1/2}(1-b\bt_my_m/B_{m-1}).  
\ee
Rewriting (\ref{St4.4}) in the form
\be\label{St4.9}
B_m = B_{m-1}(1+by_m/B_{m-1}),  
\ee
and using the inequality
$$
(1+x)^\alpha \le 1+\alpha x, \quad 0 \le \alpha \le 1, \quad x\ge 0,
$$
we get from (\ref{St4.9}) that
\be\label{St4.10}
B_m^{\bt_m} \le B_{m-1}^{\bt_m}(1+b\bt_my_m/B_{m-1}),  
\ee
 Multiplying (\ref{St4.8}) and (\ref{St4.10}) we obtain
$$
a_m^{1/2}B_m^{\bt_m} \le a_{m-1}^{1/2}B_{m-1}^{\bt_m}.
$$
Next, $B_{m-1} \ge 1$ and $t_m \le t_{m-1}$. Therefore
$$
B_{m-1}^{\bt_m} \le B_{m-1}^{\bt_{m-1}}
$$
and
\be\label{St4.11}
a_m^{1/2}b_{m}^{\bt_m} \le a_{m-1}^{1/2}B_{m-1}^{\bt_{m-1}}\le \dots
 \le a_0^{1/2}B_0^{\bt_1}  .
 \ee
Combining (\ref{St4.6}) and (\ref{St4.11}) we obtain
$$
a_m^{1+\bt_m} = \left(a_m^{1/2}B_m^{\bt_m}\right)^2 \left(a_mB_m^{-2}\right)^{\bt_m}
$$
$$
  \le \|f\|^{2}B_0^{2 \bt_1}\left((\|f\|/B_0)^{-2}+ b(2-b)\sum_{k=1}^m (ht_k)^2\right)^{-\bt_m},
$$
what completes the proof. 
\end{proof}

\section{Discussion}
\label{D}

We begin with a very simple comment on approximation by a linear method. Suppose 
$\cS := \{S_k\}_{k=1}^\infty$ is a sequence of bounded linear operators on $H$ with 
$\|S_k\| \le K$, $k=1,2,\dots$ . Then for any pair $f$, $f^\e$ satisfying $f-f^\e\|\le \e$ we have for any $k$
\be\label{D1}
\|S_k(f)-S_k(f^\e)\| \le K\e,\qquad \|f-S_k(f)\| \le \|f^\e-S_k(f^\e)\| + (K+1)\e.
\ee
Thus, the linear method $\cS$ approximates the noisy version $f$ of $f^\e$ as well as 
it approximates $f^\e$ itself up to the additive term $(K+1)\e$, which is of the order of noise. 

We now give a very simple example, which demonstrates that the PGA is not stable in the wide sense. Let $H$ be the Euclidean space $\ell^2_2$. Take the dictionary 
$\cD :=\{\bfe_1,\bfe_2\}$, $\bfe_1:=(1,0)$, $\bfe_2:=(0,1)$. Let $\e$ be a small positive number.
Consider two elements $\bv_1:= (1+\e)\bfe_1 + \bfe_2$ and $\bv_2:=  \bfe_1 + (1+\e)\bfe_2$.
Then at the first iteration the PGA will provide
$$
G_1(\bv_1) = (1+\e)\bfe_1,\quad (\bv_1)_1:=\bv_1 - G_1(\bv_1) = \bfe_2;
$$
and
$$
G_1(\bv_2) = (1+\e)\bfe_2,\quad  (\bv_2)_1:=\bv_2 - G_1(\bv_2) = \bfe_1;
$$
So, we have
$$
\|(\bv_1)_1-(\bv_2)_1\| = \sqrt{2},\qquad \|\bv_1-\bv_2\| = \e\sqrt{2}.
$$

However, results of this paper show that the PGA has stability in a special more narrow sense. Namely, if we assume that our clean signal $f^\e$ has some good properties expressed in the form $f^\e/B \in A_1(\cD)$ then the PGA approximates well the noisy version $f$ of the $f^\e$. Let us discuss this in detail and after that compare our new results 
with known results on the Weak Orthogonal Greedy Algorithm (WOGA). 
For simplicity, we discuss the case $b=1$ and $\tau =\{1\}$ ($t_k=1$, $k=1,2,...$). Then Theorem \ref{InT1} gives the following rate of convergence for the $A_1(\cD)$
\be\label{D2}
e_m(1,1) = (1+m)^{-1/6}. 
\ee
Note, that this rate of convergence result was proved in \cite{DT1}. 
Thus, in order to achieve the error $\e$ by applying the PGA to the clean signal $f^\e$ it is sufficient to make $m_\e \asymp (1/\e)^6$ iterations ($\asymp$ means "of the order"). 
Theorem \ref{T4.2} and Corollary \ref{StC1} only work with the number of iterations $m\le \e^{-2}$. In the case $m\asymp \e^{-2}$ the (\ref{D2}) gives $e_m(1,1)\asymp \e^{1/3}$. 
Corollary \ref{StC1} with $\bt = \bt(1,1,h) = h/2$ guarantees that by applying the PGA to the noisy version $f$ of $f^\e$ we achieve the following error   after (see (\ref{er1})) $m\in [(2\e)^{-2}, \e^{-2}]$ iterations
\be\label{D3}
\|f^{t,b}_m\| \le C_1\e^{\frac{\bt}{1+\bt}} = C_1\e^{\frac{h}{2+h}} .
\ee
Clearly, $\frac{h}{2+h}$ is smaller than $1/3$, but
$$
\lim_{h\to 1} \frac{h}{2+h}  = 1/3
$$
 and with $h$ close to $1$ the $\frac{h}{2+h}$  is close to $1/3$.
 
 We now proceed to the Weak Orthogonal Greedy Algorithm (WOGA). Let a sequence $\tau = \{t_k\}_{k=1}^\infty$, $0\le t_k \le 1$, be given. 
The following greedy algorithm was
defined in \cite{VT75}. We give the definition for an arbitrary dictionary $\cD$. 

 {\bf Weak Orthogonal Greedy Algorithm (WOGA($\tau$)).} Let $f_0$ be given. Then for each $m\ge 1$ we inductively define:

(1) $\varphi_m  \in \cD$ is any element satisfying
$$
|\langle f_{m-1},\varphi_m\rangle | \ge t_m
\sup_{g\in \cD} |\langle f_{m-1},g\rangle |.
$$

(2) Let $H_m := \sp (\varphi_1,\dots,\varphi_m)$ and let
$P_{H_m}(\cdot)$ denote an operator of orthogonal projection onto $H_m$.
Define
$$
G_m(f_0,\cD) := P_{H_m}(f_0).
$$

(3) Define the residual after $m$th iteration of the algorithm
$$
f_m := f_0-G_m(f_0,\cD).
$$

In the case $t_k=1$, $k=1,2,\dots$,   WOGA is called the Orthogonal
Greedy Algorithm (OGA). We now formulate the version of Theorem \ref{ST1} for the OGA. 
First of all, we point out that the WOGA is the version of the Weak Chebyshev Greedy Algorithm (WCGA) adjusted for a Hilbert space. Second, for the modulus of smoothness of a Hilbert space we have  $\rho(u)\le   u^2/2$. Then the corresponding version of Theorem \ref{ST1} reads as follows.

\begin{Theorem}\label{DiT1} Let $H$ be a  Hilbert space.   Assume that $f$ and $f^\e$ satisfy (\ref{S1}).
Then, for the OGA applied to $f$ we have
$$
\|f_m\| \le  \max\left\{2\e,\,  4(B+\e) (1+ m)^{-1/2}\right\}.
$$
\end{Theorem}

Known results (see \cite{VT75} and \cite{VTbook}) give the following rate of convergence of the OGA for the $A_1(\cD)$
\be\label{D4}
\|f_m\| \le m^{-1/2}. 
\ee
 Thus, in order to achieve the error $\e$ by applying the OGA to the clean signal $f^\e$ it is sufficient to make $m^o_\e \asymp (1/\e)^2$ iterations, which is much smaller than $m_\e\asymp (1/\e)^6$ for the PGA. 
 
 Theorem \ref{DiT1}  guarantees that by applying the OGA to the noisy version $f$ of $f^\e$ we achieve the error of order $\e$ after  
\be\label{D5}
m^o(\e) \asymp \left(\frac{1}{\e}\right)^2
\ee
iterations.
 It is an ideal situation. We note that the big difference between the behavior of $m^o(\e)$ for the OGA and of $m_\e$ for the PGA is due to very different accuracy guarantees for these two algorithms in the case of approximation of the $A_1(\cD)$. The reader can find a discussion of the rate of approximation of the PGA on the $A_1(\cD)$ in the beginning of Section 16 of \cite{VT211}. 
 
 The above discussion shows that both the PGA and the OGA, which are not stable in the wide sense (contrary to the linear methods), still "feel" the underlying clean signal $f^\e$, when we apply them to the noisy version $f$ of $f^\e$ in the case $f^\e/B\in A_1(\cD)$.

  \Addresses
 
\end{document}